\documentstyle[12pt,amsfonts,amssymb]{article}
\begin{document}
\newtheorem{theorem}{Theorem}[section]
\newtheorem{lemma}[theorem]{Lemma}
\newtheorem{defin}[theorem]{Definition}
\newtheorem{cor}[theorem]{Corollary}
\newtheorem{prop}[theorem]{Proposition}
\newtheorem{claim}[theorem]{Claim}
\newtheorem{remark}[theorem]{Remark}
\newenvironment{proof}{{\raggedright{\em Proof.}}}{\hspace{.1in}$\square$}
\newenvironment{poc}[1]{{\raggedright{\bf Proof of Claim {#1}.}}}{\hspace{
  .1in}$\square$}
\newenvironment{pol}[1]{{\raggedright{\bf Proof of Lemma {#1}.}}}{\hspace{
  .1in}$\square$}
\newenvironment{pop}[1]{{\raggedright{\bf Proof of Proposition {#1}.}}}
 {\hspace{ .1in}$\square$}
\newenvironment{pot}[1]{{\raggedright{\bf Proof of Theorem {#1}.}}}{\hspace{
  .1in}$\square$}
\newcounter{re}[theorem]
\newcounter{rom}
\newcommand{\rmn}{\stepcounter{rom}{\rm (\roman{rom})}}
\newcounter{gen}
\newcommand{\trm}[1]{\setcounter{gen}{#1}{\rm\roman{gen}}}
\newcommand{\ca}[1]{\mbox{$\Bbb{#1}$}}
\newcommand{\sca}[1]{\mbox{\scriptsize $\Bbb{#1}$}}
\newcommand{\tca}[1]{\mbox{\tiny $\Bbb{#1}$}}
\newcommand{\bol}[1]{\mbox{\boldmath ${#1}$}}
\newcommand{\Romn}[1]{\setcounter{gen}{#1}\Roman{gen}}

\title{A Height Inequality}
\author{Yuhan Zha
\\{Institute of Mathematics}\\
{Chinese Academy of Sciences} \\{yhzha1@sina.com}} \maketitle

\begin{abstract}
\label{abstract} We give a mathematical structure on an arithmetic
surface $X$, that has algebraic meanings over finite places and
can estimate the canonical norm for a relative differential form
on $X$ over infinite places. This gives a lower bound for the
canonical norm of a relative differential form on $X$, which
proves a Height Inequality on arithmetic surfaces, that implies
$abc$ conjecture.
\end{abstract}

\section{Introduction}
Let $R$ be the ring of integers of a number field $F$, and $Y={\rm
Spec}R$. Let $X$ be a stable family of curves of genus $g>1$ over
$Y$. Let \begin{equation}d_Y=\frac{1}{[F:\ca{Q}]}{\rm
log}\ca{D}_{F/\sca{Q}}\end{equation} where $\ca{D}_{F/\sca{Q}}$ is
the absolute value of the discriminant of the number field $F$.
Let $\omega_{X/Y}$ be the canonical dualizing sheaf [L1] endowed
with the canonical Hermitian metric. For a Weil divisor $D$ on
$X$, let ${\cal O}_X(D)$ be the canonical line bundle associated
with $D$ endowed with the canonical Hermitian metric [L1] on $X$.
Let $\omega_{X/Y}\cdot D$ be the product of the first arithmetic
Chern classes of $\omega_{X/Y}$ and ${\cal O}_X( D)$.

We will prove the following Theorem:
\begin{theorem}\label{ts4200}
Let $E_P$ be a curve on $X$ that is rational over $Y$. Then we
have
\begin{equation}\label{ts4210}
\omega_{X/Y}\cdot E_P\leq [F:\ca{Q}]\cdot d_Y+O(1)
\end{equation} where the constant implicit in $O(1)$ is determined by
$X_{\sca{C}}$.
\end{theorem}
This Theorem implies $abc$ conjecture and a number of other
diophantine conjectures [V1][V2].

Since our proof works for number field case, and does not work for
function field case, so we first state a difference between number
fields and function fields. Then we use this difference to
construct a mathematical structure on $X$, that has algebraic
meanings over finite places and can estimate the canonical norm of
a relative differential form on $X$ over infinite places.

The difference happens over projective lines. Let $V$ be a two
dimensional vector space over $\ca{C}$ endowed with a Hermitian
metric. Let $\ca{P}_{\sca{C}}^1$ be the complex projective line
associated with $V$. Let ${\cal O}_{\sca{P}_{\tca{C}}^1}(1)$ be
the canonical line bundle of degree $1$ over $\ca{P}_{\sca{C}}^1$.
Let $h$ be the canonical Hermitian metric on ${\cal
O}_{\sca{P}_{\tca{C}}^1}(1)$ induced from the Hermitian metric on
$V$. Let $\{x_1,\cdots, x_n\}$ be a set of different points on
$\ca{P}_{\sca{C}}^1$, where $n>2$. Then for all $1\leq i\leq n$,
there exists an element $w_{x_i}\in V$, such that
\begin{equation}\label{da8800}\left|w_{x_i}(x_i)\right|_{h}=1
\end{equation} and
\begin{equation}\label{da8810}\left|w_{x_i}(x_j)\right|_{h}<1
\end{equation} for all $j\neq i$ satisfying $1\leq j\leq n$.
Note this fact is not true over function fields of characteristic
$p>0$.

Now we use the fact above to construct the mathematical structure
mentioned above. Let $m>0$ be a positive integer.  Let $V^{\otimes
m}$ be the vector space generated by $w_1\otimes\cdots \otimes
w_m$ for all $w_i\in V$. Let $S^m(V)$ be the submodule of
$V^{\otimes m}$ that is invariant under the action of the
symmetric group  on $m$ symbols, i.e. invariant under all the
permutations of $\{w_1,\cdots,w_m\}$. Let $M_m$ be the $m+1$
dimensional vector space generated by all the sections of ${\cal
O}_{\sca{P}_{\tca{C}}^1}(m)$ over $\ca{P}_{\sca{C}}^1$. Let's
consider the homomorphism $V^{\otimes 2m} \longrightarrow M_m
\otimes M_m$ defined by mapping $w_1\otimes\cdots\otimes w_{2m}$
to $\prod_{i=1}^mw_i \otimes \prod_{i=m+1}^{2m}w_i$, where $w_i\in
V$ for all $1\leq i\leq 2m$. By restricting the map above to the
submodule $S^{2m}(V)$ of $V^{\otimes 2m}$, we have the following
map
\begin{equation}\label{da8830} S^{2m}(V)\longrightarrow
M_m \otimes M_m
\end{equation} Since the natural map $S^{2m}(V)\longrightarrow
M_{2m}$ is an isomorphism, so (\ref{da8830}) determines the
following homomorphism
\begin{equation}\label{da8840}
f_5: M_{2m}\longrightarrow M_m\otimes M_m\end{equation}

Let $\{v_0,v_1\}$ be a set of basis of $V$ over $\ca{C}$. Let
$\{\varrho_i:i\in I\}$ be a set of $n>2$ distinct complex numbers,
such that $\left|\varrho_i\right|=1$ for all $i\in I$. For $i\in
I$, let $x_i$ be a point on $\ca{P}_{\sca{C}}^1$, such that
\begin{equation}\label{da8880}\frac{v_1}{v_0}(x_i)=\varrho_i
+O(\varepsilon)\end{equation} where $\varepsilon>0$ is small
enough. Let $z$ be an analytic function on an open set $U$
containing $\{x_i:i\in I\}$ on $\ca{P}_{\sca{C}}^1$, such that
$z(x_i)=0$ for all $i\in I$, and $dz(x)\neq 0$ for all $x\in U$
satisfying $|z(x)|<r_0$, where $r_0$ is a positive number. Let
$U(r)$ be the open set in $U$ determined by $|z(x)|<r$, where
$x\in U$ and $r\in (0,r_0]$. Let $\varphi_2$ be the map from
$U(r_0)$ to the complex plane defined by $z$. Assume $\varphi_2$
is a finite morphism of degree $n$ over $|z|<r_0$.

Let $\{u_j,e_j:j\in I_1\}$ be a set of analytic functions on $U$,
such that
\begin{equation}\label{da8890}\sum_{j\in I_1}u_j\cdot e_j=0
\end{equation} Let
\begin{equation}\label{da8900}\omega= \sum_{j\in I_1}u_j\cdot de_j
\end{equation} By the definition of $f_5$, we have
\begin{equation}\label{da8910}
f_5\left(v_0^mv_1^m\right)=\sum_{l=0}^{m}b_{m,l}\cdot
v_0^{m-l}v_1^{l} \otimes v_0^{l}v_1^{m-l}\end{equation} where
$b_{m,l}\in\ca{Q}$. For $x\in U(r_0)$, assume
\begin{equation}\label{da8920} \omega=\beta_x\cdot
dz\end{equation} at point $x$, where $\beta_x\in\ca{C}$. Let $G$
be the function on $U(r_0)$ defined by the following:
\begin{equation}\label{da8930}
G=\sum_{j\in I_1}\sum_{l=0}^{m}b_{m,l}\cdot
\frac{v_1^{l}}{v_0^{l}}\cdot u_j\cdot \frac{\partial}{\partial
z}{\rm Trace}_{\varphi_2}\left(\frac{ v_1^{m-l}}{v_0^{m-l}}\cdot
e_j\right)\end{equation} Then there exists $r\in (0,r_0)$ and
$\rho\in (0,1)$, such that
\begin{equation}\label{da8940} \left|G(x)-\beta_x\cdot
\frac{v_1^m}{v_0^m}(x)\right|<\left|
\frac{v_1^m}{v_0^m}(x)\right|\cdot \rho^m
\end{equation} when $m$ is large enough. So $G$ can
estimate the norm of $\omega$ on $U(r)$. Note this fact is implied
by the existence of $w_{x_i}$ discussed in the paragraph
(\ref{da8800})-(\ref{da8810}). Next we study the algebraic meaning
of $G$ over finite places over $X$.

Let $\eta$ be the canonical section of the line bundle ${\cal
O}_X(E_P)$ vanishing along $E_P$ on $X$. Let $\xi_1$ be a section
of ${\cal O}_{X}(n_1E_P)$ on $X$, where $n_1$ is a positive
integer, such that
\begin{equation}\label{da8841}d\frac{\xi_1}{\eta^{n_1}}(x)\neq
0
\end{equation} for all the points $x\in X_{\sca{C}}$ satisfying
\begin{equation}\label{da8842}\left|\frac{\xi_1}{\eta^{n_1}}(x)
\right|<r_0\cdot \|\xi_1\|
\end{equation} where $r_0>0$ is a constant determined by
$n_1,X_{\sca{C}}$, and $\|\xi_1\|$ denotes the canonical norm of
$\xi_1$. Let $\{x_i\}$ be the set of points on $X_{\sca{C}}$,
where $\xi_1=0$. Let $\xi_2$ be a section of ${\cal
O}_{X}(n_2E_P)$ on $X$, where $n_2$ is a positive integer, such
that
\begin{equation}\label{da8843}\frac{\xi_2}{\eta^{n_2}}(x_i)
=\kappa_1\cdot \left(\varrho_i+O(\varepsilon)\right)\cdot\|\xi_2\|
\end{equation} for all $x_i$, where $\varepsilon$ is sufficiently
small, $\kappa_1$ is a positive number, and $\varrho_i\in\ca{C}$
satisfying $\left|\varrho_i\right|=1$ and $\varrho_i\neq
\varrho_j$ for $i\neq j$.

Let $\omega$ be a relative differential on $X$, such that there
exists sections $u_j$ of ${\cal O}_X(m_2E_P)$ on $X$ and sections
$e_j$ of ${\cal O}_X(m_3E_P)$ on $X$ satisfying
\begin{equation}\label{da8950} \sum_ju_j\cdot e_j=0
\end{equation}\begin{equation}\label{da8960}
\omega=\sum_j\frac{u_j}{\eta^{m_2}}\cdot d\frac{e_j}{\eta^{m_3}}
\end{equation} on $X$. In (\ref{da8930}), by replacing
$\frac{v_1}{v_0}$ with $\frac{\xi_2}{\kappa_1\cdot\|\xi_2\|\cdot
\eta^{n_2}}$, replacing $z$ with $\frac{\xi_1}{\|\xi_1\|\cdot
\eta^{n_1}}$, and replacing $u_j,e_j$ with
$\frac{u_j}{\eta^{m_2}}, \frac{e_j}{\eta^{m_3}}$ respective, we
can find a section $G_1$ of ${\cal O}_X(n_5E_P)$ on $X$ such that
\begin{equation}\label{da8970} \frac{G_1}{\eta^{n_5}}
(x)=\kappa_2\cdot G(x)
\end{equation} for all $x\in X_{\sca{C}}$ satisfying
(\ref{da8842}), where $\kappa_2$ is an element in $\ca{R}$.

Let $D_1$ be the divisor determined by $\xi_1=0$ on $X$. Let $D_3$
be the divisor determined by $G_1=0$ as a section of ${\cal
O}_X(n_5E_P)$ on $X$. Assume $D_1$ is a horizontal divisor that
does not intersect with $E_P$ on $X$, and $D_3$ is a horizontal
divisor on $X$. Let $C_2$ be the intersection cycle of $D_1$ and
$D_3$ on $X$. Let $C_3$ be the intersection cycle of $E_P$ and
$D_3$ on $X$. By considering the restriction of
$\frac{\xi_1}{\eta^{n_1}}$ to $D_3$, we see $\deg C_2$ is
determined by $\deg C_3, n_1,n_5$ and the value of
$\frac{\xi_1}{\eta^{n_1}}$ at points $x\in X_{\sca{C}}$ where
$G_1(x)=0$. Since $D_1$ does not intersect with $E_P$, so
\begin{equation}\label{da8980} \deg
C_2=\sum_i\log\left|\frac{G_1}{\eta^{n_5}}(x_i)\right|\end{equation}
where the sum is taken over all the points $x_i\in X_{\sca{C}}$
satisfying $\xi_1(x_i)=0$. So the lower bounds on $\deg C_3$ and
$\left| \frac{\xi_1}{\eta^{n_1}}(x)\right|$, where $x\in
X_{\sca{C}}$ satisfying $G_1(x)=0$, will give a lower bound on
$\sum_i\log\left|\frac{G_1}{\eta^{n_5}}(x_i)\right|$. And this
will give a lower bound on the canonical norm of $\omega$ on
$X_{\sca{C}}$.

To give a proof of {\bf Theorem~\ref{ts4200}}, we assume {\bf
Theorem~\ref{ts4200}} is not true. From this assumption we will
construct $\xi_1,\xi_2$ with similar properties discussed above
and more technical properties stated in {\bf Lemma \ref{da6999}},
such that $\deg C_3,\kappa_2$ and $\frac{\xi_1}{\eta^{n_1}}(x)$,
where $G_1(x)=0$, can be calculated. From these calculations, we
will get a lower bound for the canonical norm of $\omega$ on
$X_{\sca{C}}$. And this will lead to a contradiction.

\section{The Proof of Theorem \ref{ts4200}}
Firstly we want to prove function $G$ constructed in the
introduction has the property (\ref{da8940}).

\begin{lemma}\label{da3098}
Let $V,v_0,v_1,\ca{P}_{\sca{C}}^1,M_m,f_5,b_{m,l}$ be the elements
defined in the introduction.  Let $\{\varrho_i:i\in I\}$ be a set
of distinct complex numbers satisfying $|\varrho_i|=1$ for all
$i\in I$. Let $x_i\in\ca{P}_{\sca{C}}^1$ be the closed point
determined by
\begin{equation}\label{da0021}\frac{v_1}{v_0}(x_i)=\varrho_i
\end{equation} For $r>0$, let $U_{x_i}(r)\subset \ca{P}_{\sca{C}}^1$
be the open set determined by
\begin{equation}\label{da0022}\left|\frac{v_1}{v_0}-\varrho_i
\right|<r\end{equation} Let $m$ be a positive integer. For $x\in
\ca{P}_{\sca{C}}^1$, where $v_0(x)\neq 0$, let
$f_{4,x}\left(v_0^{m}v_1^{m}\right)$ be the rational function on
$\ca{P}_{\sca{C}}^1$ defined by the following:
\begin{equation}\label{da7021}
f_{4,x}\left(v_0^{m}v_1^{m}\right)=\sum_{l=0}^{m}b_{m,l}\cdot
\frac{v_1^{l}}{v_0^l}(x)\cdot \frac{v_1^{m-l}}{v_0^{m-l}}
\end{equation}

Then there exists $r_1,\rho_2\in (0,1)$ that satisfies the
following:
\begin{enumerate}\setcounter{rom}{0}\item[\rmn]
For all $i,j\in I$ satisfying $i\neq j$, and for all $x\in
U_{x_i}(r_1)$ and $x'\in U_{x_j}(r_1)$, we have
\begin{equation}\label{da7022}
\left|f_{4,x}\left(v_0^{m}v_1^{m}\right)(x')\right|<\rho_2^m\cdot
\left|\frac{v_1^m}{v_0^m}(x)\right|
\end{equation} for all $m>0$.
\item[\rmn] For all $i\in I$ and $x\in U_{x_i}(r_1)$ and for all
$m>0$, we have
\begin{equation}\label{da7024}
f_{4,x}\left(v_0^{m}v_1^{m}\right)(x')=\frac{v_1^m}{v_0^m}(x)
+\left(\frac{v_1}{v_0}(x')-\frac{v_1}{v_0}(x)\right)\cdot u(x')
\end{equation} for all $x'\in U_{x_i}(r_1)$,
where $u$ is an analytic function on $U_{x_i}(r_1)$.
\end{enumerate}
\end{lemma}

\begin{proof}
Let $x\in\ca{P}_{\sca{C}}^1$ be a closed point. Assume
\begin{equation}\label{da1000}\frac{v_1}{v_0}(x)=\varrho
\end{equation} where $\varrho\in \ca{C}$.
Let $v_{x,1}$ and $v_{x,0}$ be the elements in $V$, such that
\begin{equation}\label{da8342}v_{x,1}=\frac{1}{\left(1+\left|
\varrho\right|^2\right)^{\frac{1}{2}}}\cdot\left( v_1-\varrho\cdot
v_0\right)
\end{equation}
\begin{equation}\label{da8343}v_{x,0}=
\frac{1}{\left(1+\left|\varrho
\right|^2\right)^{\frac{1}{2}}}\cdot\left( \overline{\varrho}\cdot
v_1+v_0\right)
\end{equation} on $\ca{P}_{\sca{C}}^1$.
By (\ref{da8342}) (\ref{da8343}), we have
\begin{equation}v_0=\frac{1}{
\left(1+\left|\varrho\right|^2\right)^{\frac{1}{2}}}
\cdot\left(v_{x,0}-\overline{\varrho}\cdot v_{x,1}\right)
\label{da9680}\end{equation}
\begin{equation}\label{da9690}v_1=\frac{1}{
\left(1+\left|\varrho\right|^2\right)^{\frac{1}{2}}}
\cdot\left(\varrho\cdot v_{x,0}+v_{x,1}\right)
\end{equation}

Note \begin{equation}\label{da9691}
\left|\frac{v_{x,0}}{v_0}\right|^2+\left|\frac{v_{x,1}}{v_0}
\right|^2=\left|\frac{v_1}{v_0}\right|^2+1\end{equation} on
$\ca{P}_{\sca{C}}^1$. So  for $l=0,1$ and $x'\in U_{x_j}(r_1)$, we
have
\begin{equation}\label{da0060}\left|\frac{v_{x,l}}{v_0}(x')
\right|\leq\left(1+\left|\frac{v_1}{v_0}(x')
\right|^2\right)^{\frac{1}{2}}
\end{equation}

Note for $0\leq l\leq {m}$, we have
\begin{eqnarray}\lefteqn{ f_5\left(
\frac{(2m)!}{(2m-l)!\cdot l!}\cdot
v_{x,0}^{2m-l}v_{x,1}^{l}\right)=}\nonumber\\
&&\sum_{l_1=0}^{l} \frac{m!}{(m-l_1)!\cdot l_1!}\cdot
\frac{m!}{(m-l+l_1)!\cdot
(l-l_1)!}\cdot\nonumber\\
&&\frac{v_{x,0}^{m-l+l_1} v_{x,1}^{l-l_1}}{v_0^{m}}\otimes
\frac{v_{x,0}^{m-l_1} v_{x,1}^{l_1}}{v_0^{m}}
\label{da9801}\end{eqnarray} And for ${m}<l\leq 2m$, we have
\begin{eqnarray}\lefteqn{ f_5\left(
\frac{(2m)!}{(2m-l)!\cdot l!}\cdot
v_{x,0}^{2m-l}v_{x,1}^{l}\right)=}\nonumber\\
&&\sum_{l_1=l}^{2m} \frac{m!}{(2m-l_1)!\cdot (l_1-m)!}\cdot
\frac{m!}{(l_1-l)!\cdot (m+l-l_1)!}
\cdot\nonumber\\
&& \frac{v_{x,0}^{2m-l_1}
v_{x,1}^{l_1-m}}{v_0^{m}}\otimes\frac{v_{x,0}^{l_1-l}
v_{x,1}^{m+l-l_1}}{v_0^{m}}\label{da3340}\end{eqnarray}

Let $f_{4,x}$ be the linear map from $M_{2m}$ to the vector space
of rational functions on $\ca{P}_{\sca{C}}^1$ defined by $f_5$
followed by the map from $M_m\otimes M_m$ to rational functions on
$\ca{P}_{\sca{C}}^1$ that maps $t_1\otimes t_2$ to
$\frac{t_1}{v_0^m}(x)\cdot \frac{t_2}{v_0^m}$. Since
$v_{x,1}(x)=0$ and
$\frac{v_{x,0}}{v_0}(x)=\left(1+|\varrho|^2\right)^{\frac{1}{2}}$,
so by (\ref{da9801}) (\ref{da3340}), we have
\begin{eqnarray}
f_{4,x}\left( v_{x,0}^{2m-l}v_{x,1}^{l}\right)= \frac{(2m-l)!\cdot
l!}{(2m)!}\cdot\left(1+|\varrho|^2\right)^{\frac{m}{2}}\cdot
\frac{m!}{(m-l)!\cdot l!}\cdot
\frac{v_{x,0}^{m-l}v_{x,1}^l}{v_0^m}\label{da0100}\end{eqnarray}
when $0\leq l\leq m$, and $f_{4,x}\left(
v_{x,0}^{2m-l}v_{x,1}^{l}\right)$ is equal to $0$ when $l>m$. Note
\begin{equation}\label{da0101}
\frac{(2m-l)!\cdot l!}{(2m)!}\cdot\frac{m!}{(m-l)!\cdot l!} \leq
2^{-l}\end{equation} for $0\leq l\leq m$, and
\begin{equation}\label{da0103}
\frac{(2m-l)!\cdot l!}{(2m)!}\cdot\frac{m!}{(m-l)!\cdot l!} \leq
2^{-\frac{m}{2}}\cdot
3^{-l+\frac{m}{2}}=\frac{3^{\frac{m}{2}}}{2^{\frac{m}{2}}}\cdot
3^{-l}\end{equation} for $\left[\frac{m}{2}\right]+1\leq l\leq m$.

Let $r_1$ be a positive number. Assume $x$ is a point in
$U_{x_i}(r_1)$. Then we have
\begin{eqnarray}\lefteqn{
v_0^{{m}}v_1^{{m}}= \left(1+\left|\varrho\right|^2\right)^{-{m}}
\cdot\left(v_{x,0}-\overline{\varrho}\cdot
v_{x,1}\right)^m\cdot\left(\varrho\cdot v_{x,0}+v_{x,1}\right)^m
}\nonumber\\
&=&\left(1+\left|\varrho\right|^2\right)^{-{m}}
\cdot\left(\varrho\cdot
v_{x,0}^2+\left(1-\left|\varrho\right|^2\right)\cdot
v_{x,0}v_{x,1}-\overline{\varrho}\cdot v_{x,1}^2
\right)^m\label{da3300}\end{eqnarray} Let
\begin{equation}\label{da3310}\lambda_1=\left|\frac{v_{x,1}}{
v_{x,0}}\right|
\end{equation} Consider the subset determined by
 $\lambda_1\leq {2}$ over $U_{x_j}(r_1)$, where
$j\neq i$. By (\ref{da0100}) (\ref{da0101}), we have
\begin{eqnarray}\lefteqn{\left|f_{4,x}(v_{x,0}^{2m-l}v_{x,1}^{l}
)\right| \leq 2^{-l}\cdot \left(1+|\varrho|^2\right)^{\frac{m}{2}}
\cdot\left|\frac{v_{x,0}^{m-l}v_{x,1}^l}{v_0^m}\right|}\nonumber\\
&\leq& 2^{-l}\cdot \left(1+|\varrho|^2\right)^{\frac{m}{2}}
\cdot\left|\frac{v_{x,0}^{m}}{v_0^m}\right|\cdot\lambda_1^l
\hspace{1in}\label{da3320}\end{eqnarray} By (\ref{da9691}), we
have
\begin{equation}\label{da3330} \left|\frac{v_{x,0}}{v_0}\right|^2
\cdot \left(1+\lambda_1^2\right)=\left|\frac{v_1}{v_0}\right|^2+1
\end{equation}
By $\lambda_1\leq {2}$, we have
\begin{equation}\label{da3350}
\left(1+\frac{\lambda_1^2}{4}\right)^2<1+\lambda_1^2
\end{equation} By $j\neq i$, we have $\lambda_1>0$.
So when $r_1>0$ is small enough, we have
\begin{eqnarray}\lefteqn{\left|\frac{v_{x,0}^m}{v_0^m}\right|
\cdot
\left(\left|\varrho\right|+\frac{\left|1-\left|\varrho\right|^2
\right|\cdot\lambda_1}{2}+\frac{\left|\varrho\right|
\cdot\lambda^2_1}{4}\right)^m}\nonumber\\
&<&\left|\varrho\right|^m \cdot
\left(1+\lambda_1^2\right)^\frac{m}{2} \cdot
\left|\frac{v_{x,0}^m}{v_0^m}\right|\cdot \rho_1^m
\nonumber\\
&=&\left|\varrho\right|^m \cdot
\left(\left|\frac{v_1}{v_0}\right|^2+1\right)^{\frac{m }{2}}\cdot
\rho_1^m\label{da3360}\end{eqnarray} where $\rho_1\in (0,1)$. So
when we take $r_1$ small enough, such that
\begin{equation}\label{da3370}\left(1+(1-r_1)^2\right)^{-1}
\cdot \left(1+(1+r_1)^2\right)\cdot\rho_1^2<1\end{equation} by
(\ref{da3300}) (\ref{da3320}) (\ref{da3360}) (\ref{da3370}), we
have
\begin{eqnarray}\lefteqn{\left|f_{4,x}(v_0^{{m}}v_1^{{m}})\right|
\leq \left(1+\left|\varrho\right|^2\right)^{-\frac{m}{2}}
\cdot\left|\frac{v_{x,0}^m}{v_0^m}\right| }\nonumber\\
&&\cdot
\left(\left|\varrho\right|+\frac{\left|1-\left|\varrho\right|^2
\right|\cdot\lambda_1}{2}+\frac{\left|\varrho\right|
\cdot\lambda_1^2}{4}\right)^m\nonumber\\
&<&\left(1+\left|\varrho\right|^2\right)^{-\frac{m}{2}}
\cdot\left|\varrho\right|^m
\cdot \left(\left|\frac{v_1}{v_0}\right|^2+1\right)^{\frac{m }{2}}
\cdot \rho_1^m\nonumber\\
&<&\rho_2^m \cdot\left|\varrho\right|^{{m}}
\label{da0110}\end{eqnarray} over $U_{x_j}(r_1)$, where $j\neq i$
and $\rho_2\in (0,1)$.

Now consider the subset determined by
 $\lambda_1\geq {2}$ over $U_{x_j}(r_1)$, where
$j\neq i$. By (\ref{da9691}), we have
\begin{equation}\label{da3450}\left|
\frac{v_{x,1}}{v_0}\right|^2\cdot \left(1+\lambda_1^{-2}\right)
=1+\left|\frac{v_1}{v_0}\right|^2\end{equation} So we have
\begin{equation}\label{da3460}
\left| \frac{v_{x,1}}{v_0}\right|^2\geq \frac{4}{5} \cdot
\left(1+\left|\frac{v_1}{v_0}\right|^2\right)\end{equation}
Therefore over $U_{x_j}(r_1)$, by (\ref{da9691}), we have
\begin{equation}\label{da3470}
\left|\frac{v_{x,0}}{v_0}\right|^2\leq \frac{1}{5} \cdot
\left(1+\left|\frac{v_1}{v_0}\right|^2\right) < \frac{1}{5} \cdot
\left(1+(1+r_1)^2\right)\end{equation}

Assume \begin{equation}\label{da3490} v_0^mv_1^m =\sum_{i=0}^{2m}
b_{m,i,x}\cdot v_{x,0}^{2m-i}v_{x,1}^i\end{equation} where
$b_{m,i,x}\in\ca{C}$. By (\ref{da0100}) (\ref{da0103})
(\ref{da3470}) (\ref{da0060}), over $U_{x_j}(r_1)$ for
$\left[\frac{m}{2}\right]+1\leq l\leq 2m$, we have
\begin{eqnarray}\lefteqn{
\left|f_{4,x}\left( v_{x,0}^{2m-l}v_{x,1}^{l}\right)\right| <
}\nonumber\\
&&\frac{3^\frac{m}{2}}{2^\frac{m}{2}}\cdot 3^{-l}\cdot
(1+(1+r_1)^2)^{\frac{m}{2}}
\cdot\frac{1}{5^\frac{m-l}{2}}\cdot (1+(1+r_1)^2)^{\frac{m}
{2}}\nonumber\\
&=&\frac{3^\frac{m}{2}}{10^\frac{m}{2}}\cdot
(1+(1+r_1)^2)^{{m}}\cdot
\left(\frac{5}{9}\right)^\frac{l}{2}\label{da3480}\end{eqnarray}
By (\ref{da3490}) (\ref{da3480}) (\ref{da3300}), we have
\begin{eqnarray}\lefteqn{
\sum_{i=\left[\frac{m}{2}\right]+1}^{2m}
\left|b_{m,i,x}\cdot f_{4,x}\left(v_{x,0}^{2m-i}v_{x,1}^i
\right)\right| <}\nonumber\\
&&\frac{\left(1+(1+r_1)^2\right)^{{m}}
}{\left(1+\left|\varrho\right|^2\right)^{{m}}} \cdot
\frac{3^\frac{m}{2}}{10^\frac{m}{2}}\cdot
\left(|\rho|+\frac{5^{\frac{1}{2}}\cdot
\left|1-|\rho|^2\right|}{3}+ \frac{5\cdot |\rho|}{9}\right)^m
\label{da3440}\end{eqnarray} When $r_1\longrightarrow 0^+$, we
have $|\rho|\longrightarrow 1$. So (\ref{da3440}) implies
\begin{equation}\label{da3500}
\sum_{i=\left[\frac{m}{2}\right]+1}^{2m}\left|b_{m,i,x}\cdot
f_{4,x}\left(v_{x,0}^{2m-i}v_{x,1}^i \right)\right|
<|\rho|^m\cdot\rho_4^m
\end{equation}
where $\rho_4\in (0,1)$, when $r_1>0$ is small enough. By
(\ref{da3470}) (\ref{da0060}) (\ref{da0100}) (\ref{da0101}), over
$U_{x_j}(r_1)$ for $0\leq l\leq \left[\frac{m}{2}\right]$, we have
\begin{equation}\label{da3510}
\left|f_{4,x}\left( v_{x,0}^{2m-l}v_{x,1}^{l}\right)\right| \leq
\frac{1}{5^\frac{m}{4}}\cdot 2^{-l}\cdot (1+(1+r_1)^2)^{{m}}
\end{equation} By (\ref{da3510}) (\ref{da3300}), we have
\begin{eqnarray}\lefteqn{\sum_{i=0}^{\left[\frac{m}{2}\right]}
\left|b_{m,i,x}\cdot f_{4,x}\left(v_{x,0}^{2m-i}v_{x,1}^i
\right)\right| <}\nonumber\\
&&\frac{\left(1+(1+r_1)^2\right)^{{m}}
}{\left(1+\left|\varrho\right|^2\right)^{{m}}} \cdot
\frac{1}{5^\frac{m}{4}}\cdot \left(|\rho|+
\frac{\left|1-|\rho|^2\right|}{2}+ \frac{|\rho|}{4}\right)^m
\label{da3520}\end{eqnarray} When $r_1\longrightarrow 0^+$, we
have $|\rho|\longrightarrow 1$. So (\ref{da3520}) implies
\begin{equation}\label{da3530}
\sum_{i=0}^{\left[\frac{m}{2}\right]}\left|b_{m,i,x}\cdot
f_{4,x}\left(v_{x,0}^{2m-i}v_{x,1}^i \right)\right|
<|\rho|^m\cdot\rho_5^m
\end{equation}
where $\rho_5\in (0,1)$, when $r_1>0$ is small enough. Then
(\ref{da0110}) (\ref{da3500}) (\ref{da3530}) imply (\trm{1}) is
true.

Over $U_{x_i}(r_1)$, by (\ref{da3300}), we have
\begin{equation}
f_{4,x}(v_0^{{m}}v_1^{{m}})=\frac{\varrho^m\cdot
f_{4,x}(v_{x,0}^{2m})}{ \left(1+\left|\varrho\right|^2\right)^{m}}
+\sum_{i=1}^{2m}b_{i,2}\cdot f_{4,x}(v_{x,0}^{2m-i}v_{x,1}^i)=
\varrho^{{m}} +\left(\frac{v_1}{v_0}-\varrho\right)\cdot u
\label{da0120}\end{equation} where $b_{i,2}\in\ca{C}$ and $u$ is
an analytic function on $U_{x_i}(r_1)$. So (\trm{2}) is true.
\end{proof}

\begin{theorem}\label{da9000} Let $\{x_i:i\in I\},
U_{x_i}(r),r_1,b_{m,l},v_0,v_1$ be the elements defined in  Lemma
{\bf \ref{da3098}}. Let $r_2>0$ be a constant. Let $z$ be an
analytic function on $\bigcup_{i\in I}U_{x_i}(r_1)$ in
$\ca{P}_{\sca{C}}^1$ that satisfies the following:
\begin{enumerate}\setcounter{rom}{0}
\item[\rmn] $|z(x_i)|<r_2$ for all $i\in I$. \item[\rmn]
$dz(x)\neq 0$ for all $x\in \bigcup_{i\in I} U_{x_i}(r_1)$
satisfying $\left|z(x)\right|<r_2$. \item[\rmn] For all $i\in I$
and $\alpha\in\ca{C}$ satisfying $|\alpha|<r_2$, then there exists
a unique point $x\in U_{x_i}(r_1)$, such that $z(x)=\alpha$.
\end{enumerate}

Let $\{u_i,e_i:i\in I_1\}$ be a set of analytic functions on $
\bigcup_{i\in I}U_{x_i}(r_1)$ satisfying
\begin{equation}\label{da9010}\sum_{i\in I_1}u_i\cdot e_i=0
\end{equation} Let \begin{equation}\label{da9020}
\omega=\sum_{i\in I_1}u_i\cdot de_i\end{equation}  Let $U'\subset
\bigcup_{i\in I}U_{x_i}(r_1)$ be the subset determined by
\begin{equation}\label{da9050} |z|<r_2\end{equation} For $x\in
U'$, assume
\begin{equation}\label{da9030}\omega= \beta_x\cdot dz
\end{equation} at point $x$, where $\beta_x\in\ca{C}$.
Let $\varphi_2$ be the map from $U'$ to the
complex plane defined by $z$. Let
\begin{equation}\label{da9040}
G=\sum_{j\in I_1}\sum_{l=0}^{m}b_{m,l}\cdot
\frac{v_1^{l}}{v_0^{l}}\cdot u_j\cdot \frac{\partial}{\partial
z}{\rm Trace}_{\varphi_2}\left(\frac{ v_1^{m-l}}{v_0^{m-l}}\cdot
e_j\right)\end{equation} on $U'$. Assume $\beta_x\neq 0$ for all
$x\in U'$. Then there exists $\rho_6\in (0,1)$ and $N_5>0$, such
that
\begin{equation}\label{da9060} \left|G(x)-\beta_x\cdot
\frac{v_1^m}{v_0^m}(x)\right|<\left|
\frac{v_1^m}{v_0^m}(x)\right|\cdot\rho_6^m
\end{equation} for all $x\in U'$ and $m>N_5$.
\end{theorem}

\begin{proof} For $x\in U'$, there exists a set of analytic
functions  $\{e_{x,0},e_{x,1},\cdots,e_{x,n_3} \}$ on
$\bigcup_{i\in I}U_{x_i}(r_1)$, that satisfies the following:
\begin{enumerate}\item ${e_{x,0}}(x)=1$. \item
${e_{x,1}}(x)=0$, and $d{e_{x,1}}-dz$ vanishes at point $x$, \item
$e_{x,i}$ vanishes at point $x$ with order $2$ for all $2\leq
i\leq n_3$. \item There exists a set of analytic functions
$\{u_{x,0},\cdots,u_{x,n_3}\}$  on
$\bigcup_{i\in I}U_{x_i}(r_1)$,
such that
\begin{equation}\label{da8200}\sum_{i\in I_1}u_i\otimes e_i
=\sum_{i=0}^{n_3}u_{x,i}\otimes e_{x,i}
\end{equation}
\end{enumerate}

Let \begin{equation}G_x= \sum_{i\in I_1}\sum_{l=0}^{m}b_{m,l}\cdot
\frac{v_1^{l}}{v_0^{l}}(x)\cdot {u_i}(x) \cdot
\frac{v_1^{m-l}}{v_0^{m-l}} \cdot {e_i}
\label{da8131}\end{equation} Since
\begin{equation}\label{da8210}\sum_{i=0}^{n_3}u_{x,i}\cdot
e_{x,i}=0
\end{equation} and $e_{x,i}(x)=0$ for $i>0$, so
\begin{equation}\label{da8220} u_{x,0}(x)=0
\end{equation} By (\ref{da8220}), we have
\begin{eqnarray}\lefteqn{
G_x= \sum_{i=1}^{n_3}\sum_{l=0}^{m}b_{m,l}\cdot
\frac{v_1^{l}}{v_0^{l}}(x)\cdot {u_{x,i}}(x) \cdot
\frac{v_1^{m-l}}{v_0^{m-l}} \cdot {e_{x,i}}
}\nonumber\\
&=&f_{4,x}\left(\frac{v_1^m}{v_0^m}\right)\cdot \sum_{i=1}^{n_3}
{u_{x,i}}(x)\cdot {e_{x,i}}
\hspace{.7in}\label{da8132}\end{eqnarray}

Since
\begin{equation}\label{da8230} \sum_{i=0}^{n_3}
{u_{x,i}}\cdot d{e_{x,i}}=\omega\end{equation} and $e_{x,i}$
vanishes at $x$ with order $2$ for $i>1$, so we have
\begin{equation}\label{da8250}{u_{x,1}}
(x)=\beta_x
\end{equation}

By {\bf Lemma \ref{da3098}}, we have
\begin{equation}\label{da8251}
\left|\frac{\partial}{\partial z}G_x\right|<
\left|\frac{v_1}{v_0}(x)\right|^m\cdot\rho_3^m
\end{equation} on $U'\cap U_{x_j}(r_1)$
when $x$ is not a point in $U_{x_j}(r_1)$, where $\rho_3\in
(0,1)$. By {\bf Lemma \ref{da3098}}, and ${e_{x,i}}$ vanishes at
$x$ with order $2$ for $i>1$, and ${e_{x,1}}$ vanishes at $x$, and
\begin{equation}\label{da8252}\frac{\partial
{e_{x,1}}}{\partial z} (x)=1
\end{equation} we have the value of $\frac{\partial}{\partial
z}\left({\rm Trace}_{\varphi_2}G_x\right)$ at point $x$ is equal
to $\frac{v_1^m}{v_0^m}(x)\cdot\left({\beta_x}+O(\rho_3^m)\right)$
for $x\in U'$. Therefore we have
\begin{equation}\label{da8253}
G(x)=\frac{v_1^m}{v_0^m}(x)\cdot\left({\beta_x}+O(\rho_3^m)\right)
\end{equation} for $x\in U'$. So our Theorem is true.
\end{proof}

\hspace{1in}

Now we want to construct $\xi_1,\xi_2$ that have the properties
stated in the introduction. Let $\eta$ be the canonical section of
${\cal O}_X(E_P)$ vanishing along $E_P$ on $X$. Let $S$ be the set
of all complex embedding of $F$ in $\ca{C}$. For $\sigma\in S$,
let $x_{0,\sigma}$ be the point determined by
$E_P\otimes_\sigma\ca{C}$ on $X_\sigma$. Let $s_{0,\sigma}$ be an
element in
\begin{equation}\label{da1220}R^0\pi_*(X,{\cal
O}_X(3gE_P))\otimes_\sigma\ca{C}\end{equation} such that the norm
of $s_{0,\sigma}$ under the canonical Hermitian metric on
(\ref{da1220}) is equal to $1$ and $s_{0,\sigma}$ is orthogonal to
the subspace of (\ref{da1220}) vanishing at $x_{0,\sigma}$. Let
$s_{1,\sigma}$ be an element in the subspace of (\ref{da1220})
vanishing at $x_{0,\sigma}$ such that the norm of $s_{1,\sigma}$
is equal to $1$ and $s_{1,\sigma}$ is orthogonal to the subspace
of (\ref{da1220}) vanishing at $x_{0,\sigma}$ with order $2$. Let
$t_{\sigma}$ be the analytic function on an open set of $X_\sigma$
containing $x_{0,\sigma}$ defined by
\begin{equation}\label{da1210}
t_\sigma=\frac{s_{1,\sigma}}{s_{0,\sigma}}
\end{equation} We define a Hermitian metric $h_1$ on $\pi_*({\cal
O}_X(nE_P)/{\cal O}_X((n-l)E_P))$, where $l>0$, such that
\begin{equation}\label{da1200}
\left<t_\sigma^{-i},{t_\sigma^{-j}}\right>_{h_1}
=\delta_{i,j}\end{equation} where $\delta_{i,j}$ is the number
that is equal $0$ when $i\neq j$ and is equal to $1$ when $i=j$.

\begin{lemma}\label{ts4780} There exists $N_0>0$ determined by
$X_{\sca{C}}$, such that when \begin{equation}\label{da1201}
\frac{1}{[F:\ca{Q}]}\cdot\omega_{X/Y}\cdot
E_P>d_Y+N_0\end{equation} then we have the following:

Let $n$ be a positive integer. Let $\zeta_n$ be a nonzero section
of the restriction of ${\cal O}_X(nE_P)$ to $E_P$, such that
\begin{equation}\label{dh1230}\log\|\zeta_n\|_{\sigma_1}\leq
\log\|\zeta_n\|_{\sigma_2}+2
\end{equation} for all $\sigma_1,\sigma_2\in S$, where $\|\cdot\|$
denotes the norm under $h_1$. Let $C_0$ be the divisor determined
by $\zeta_n=0$ as a section of ${\cal O}_X(nE_P)$ on $E_P$. The
push-forward cycle of $C_0$ on $Y$ is still denoted by $C_0$.

Then there exists an element $\zeta_{n}'$ in $\pi_*({\cal
O}_X(nE_P)/{\cal O}_X )\otimes {\cal O}_Y(-C_0)$ such that
$\zeta_n'-\zeta_n=0$ as a section of the restriction of ${\cal
O}_X(nE_P)$ to $E_P$ and \begin{equation}\label{dh1240}
\log\|\zeta_n'\|_\sigma= \log\|\zeta_n\|_\sigma+O(1)\end{equation}
for all $\sigma\in S$, where the constant implicit in $O(1)$ is
determined by $X_{\sca{C}},n$.
\end{lemma}

\begin{proof}
By induction, we assume there exists an element $\zeta_{n,i}$ in
\begin{equation}\label{ts9330}\pi_*({\cal
O}_X(nE_P)/{\cal O}_X((n-i)E_P))\otimes {\cal O}_Y(-C_0)
\end{equation} where $i\geq 1$, such that
\begin{equation}\label{dh1280}\zeta_{n,i}-
\zeta_{n,i-1}=0
\end{equation} as an element in $\pi_*({\cal
O}_X(nE_P)/{\cal O}_X((n-i+1)E_P))\otimes {\cal O}_Y(-C_0)$, and
\begin{equation}\log \|\zeta_{n,i}\|_\sigma=\log \|\zeta_{n}\|_\sigma
+O(1) \label{dh1290}\end{equation} for all $\sigma\in S$.

Note we have exact sequence
\begin{eqnarray}\lefteqn{0\longrightarrow {\cal
O}_X((n-i)E_P)/{\cal
O}_X((n-i-1)E_P)}\nonumber\\
&&\longrightarrow
{\cal O}_X(nE_P)/{\cal O}_X((n-i-1)E_P)\nonumber\\
&& \longrightarrow {\cal O}_X(nE_P)/{\cal
O}_X((n-i)E_P)\longrightarrow 0 \label{ts9332}\end{eqnarray} Since
\begin{equation}\label{ts9340} \widehat{c}_1({\cal
O}_X((n-i)E_P)|_{E_P})=-(n-i)\omega_{X/Y}\cdot E_P\end{equation}
so there exists a set of elements $\{v_{i,j}\}$ in
\begin{equation}\label{ts9350}({\cal O}_X((n-i)E_P)/{\cal
O}_X((n-i-1)E_P))\otimes{\cal O}_Y(-C_0)\end{equation} such that
$\{v_{i,j}\}$ is a set of generators of (\ref{ts9350}) as a
$\ca{Z}$ module, and
\begin{equation}\label{ts9360}\log\|v_{i,j}\|\leq
\frac{1}{[F:\ca{Q}]}((n-i)\omega_{X/Y}\cdot
E_P+\deg(C_0))+d_Y+O(1)
\end{equation}  By (\ref{dh1230}), we have
\begin{equation}\label{ts9370}\log\|\zeta_n\|_\sigma
=\frac{1}{[F:\ca{Q}]}(n\cdot\omega_{X/Y}\cdot
E_P+\deg(C_0))+d_Y+O(1)
\end{equation} By (\ref{dh1290}) (\ref{ts9370}) (\ref{ts9360})
(\ref{da1201}), we have \begin{equation}\label{ts9380} \log
\|\zeta_{n,i}\|_\sigma\gg \log\|v_{i,j}\|\end{equation} Then by
(\ref{ts9332}), there exists $\zeta_{n,i+1}$ in
\begin{equation}\label{ts9400}\pi_*({\cal
O}_X(nE_P)/{\cal O}_X((n-i-1)E_P))\otimes {\cal O}_Y(-C_0)
\end{equation} such that $\zeta_{n,i+1}-\zeta_{n,i}=0$ in
(\ref{ts9330}), and \begin{equation}\label{dh1310}
\log\|\zeta_{n,i+1}\|=\log\|\zeta_n\|+O(1)\end{equation}  Let
$\zeta_n'=\zeta_{n,n}$. We see our Lemma is true.
\end{proof}

\hspace{1in}

Consider the natural map by restriction:
\begin{equation}\label{da0830} R^0\pi_*(X,{\cal
O}_X(3gE_P))\longrightarrow \pi_*({\cal O}_X(3gE_P)/{\cal
O}_X)\end{equation} Let $W$ be the quotient of map (\ref{da0830}).
Note for all $e_1\in\pi_*({\cal O}_X(3gE_P)/{\cal O}_X)$ and
$\omega_1\in R^0\pi_*(X,\omega_{X/Y})$
\begin{equation}\label{da0852}{\rm Res}_{E_P/Y}\left(
\frac{e_1}{\eta^{3g}}\cdot \omega_1\right)
\end{equation} is an element in $R$. So we have a natural map
\begin{equation}\label{da0820}
\pi_*({\cal O}_X(3gE_P)/{\cal O}_X)\longrightarrow
R^0\pi_*(X,\omega_{X/Y})^\vee\end{equation} defined by
(\ref{da0852}), where $R^0\pi_*(X,\omega_{X/Y})^\vee$ denotes the
dual of $R^0\pi_*(X,\omega_{X/Y})$. Note map (\ref{da0820})
vanishes on $R^0\pi_*(X,{\cal O}_X(3gE_P))$. So we have a natural
map
\begin{equation}\label{da0810}
W\longrightarrow R^0\pi_*(X,\omega_{X/Y})^\vee
\end{equation} induced from (\ref{da0820}).
Let $C$ be the degeneracy cycle of (\ref{da0810}).

\begin{theorem}\label{dh0800}
Let $n_0>9g^2$ be an integer. Let $C$ be the cycle on $Y$ defined
above. Let $\zeta_{n_0}$ be a section of the restriction of ${\cal
O}_X(n_0E_P)$ to $E_P$, that satisfies
\begin{equation}\label{ts4771}\log\|\zeta_{n_0}\|_{\sigma_1} \geq
\log\|\zeta_{n_0}\|_{\sigma_2}+2\end{equation} for all
$\sigma_1,\sigma_2\in S$, where $\|\cdot\|$ denotes the norm under
$h_1$. Let $C_0$ be the divisor determined by $\zeta_{n_0}=0$ on
$E_P$. Then there exists $N_0>0$, which is determined by
$X_{\sca{C}},n_0$, that satisfies the following:

If the rational curve $E_P\subset X$ satisfies
\begin{equation}\label{ts4770}\frac{1}{[F:\ca{Q}]}\cdot
\omega_{X/Y}\cdot E_P>d_Y+N_0
\end{equation} then
there exists an element $\xi$ in $R^0\pi_*(X, {\cal O}_X(n_0E_P))
\otimes {\cal O}_Y(-C_0+C)$ that satisfies the following:
\begin{enumerate}\setcounter{rom}{0}\item[\rmn]
The restriction of $\xi$ to $E_P$ is equal to $\zeta_{n_0}$.
\item[\rmn] For all $\sigma\in S$, we have
\begin{equation}\label{dh1400}\log|\xi(x_{0,\sigma})|_h>
\log\|\xi\|+O(1)
\end{equation}  where $\|\xi\|$ denotes
the norm of $\xi$ under the canonical Hermitian metric on
$R^0\pi_*(X, {\cal O}_X(n_0E_P))$, and $h$ denotes the canonical
Hermitian metric on ${\cal O}_X(n_0E_P)$, and the constant
implicit in $O(1)$ is determined by $X_{\sca{C}},n_0$.
\end{enumerate}
\end{theorem}

\begin{proof}
By Arithmetic Riemann-Roch Theorem,
\begin{equation}\widehat{c}_1(R^0\pi_*(X,{\cal
O}_X(3gE_P)),h_2)\geq -\frac{3g(3g+1)}{2}\omega_{X/Y}\cdot
E_P+O(1) \label{dh0850}\end{equation} where $h_2$ is the canonical
Hermitian metric. Moreover we have
\begin{equation}\label{dh0851}\log\|e\|\geq O(1)
\end{equation} for all nonzero $e\in R^0\pi_*(X,{\cal
O}_X(3gE_P))$. So there exists a set of elements $\{e_i'\}$ in
\begin{equation}\label{dh0852}
R^0\pi_*(X,{\cal O}_X(3gE_P))\otimes {\cal
O}_Y(-C_0+C)\end{equation} such that $\{e_i'\}$ generate
(\ref{dh0852}) over $\ca{Z}$, and for all $e_i'$, we have
\begin{equation}\label{da0840} \log\|e_i'\|<
\frac{3g(3g+1)}{2[F:\ca{Q}]}\omega_{X/Y}\cdot E_P+ \frac{\deg
C_0}{[F:\ca{Q}]}+3g\cdot d_Y+O(1)
\end{equation}

Let
 $\zeta_{n_0,n_0}$ be the element
constructed from $\zeta_{n_0}$ in {\bf Lemma \ref{ts4780}}. Note
for all $\omega_1\in R^0\pi_*(X,\omega_{X/Y})$,
\begin{equation}\label{da0860}{\rm Res}_{E_P/Y}\left(\omega_1
\cdot \zeta_{n_0,n_0}\right)
\end{equation} is a section of ${\cal O}_Y(-C_0)$ on $Y$.
So there exists unique $w_1$ in $W\otimes {\cal O}_Y(-C_0+C)$,
such that
\begin{equation}\label{da0800} {\rm Res}_{E_P/Y}\left(\omega_1
\cdot \zeta_{n_0,n_0} \right)={\rm Res}_{E_P/Y}\left(\omega_1\cdot
w_1\right)
\end{equation} for all $\omega_1\in R^0\pi_*(X,\omega_{X/Y})$.
By (\ref{ts4771}), we have \begin{equation}\label{da0790}
\log\|\zeta_{n_0,n_0}\|_\sigma=\frac{n_0
}{[F:\ca{Q}]}\omega_{X/Y}\cdot E_P+\frac{\deg C_0}{[F:\ca{Q}]}
+O(1)\end{equation} for all $\sigma\in S$, so
\begin{equation}\label{da0780}
\log\|w_1\|_\sigma=\frac{n_0}{[F:\ca{Q}]}\omega_{X/Y}\cdot
E_P+\frac{\deg C_0}{[F:\ca{Q}]} +O(1)
\end{equation} where $\|w_1\|$ denotes the norm under natural
Hermitian metric on $W$. By (\ref{da0840}) (\ref{da0780}), we have
\begin{equation}\label{da0770} \log\|w_1\|_\sigma>
\log\|e_{i}'\|
\end{equation} for all $e_i'$. Since $\{e_i'\}$ generate
(\ref{dh0852}), so there exists $w_2\in \pi_*({\cal
O}_X(3gE_P)/{\cal O}_X)\otimes {\cal O}_Y(-C_0+C)$, such that
\begin{equation}\label{da0760} \log\|w_2\|<
\log\|\zeta_{n_0,n_0}\| +O(1)\end{equation} and
\begin{equation}\label{da0750}
{\rm Res}_{E_P/Y}\left(\omega_1 \cdot \left(\zeta_{n_0,n_0}
-w_2\right) \right)=0
\end{equation} for all $\omega_1\in R^0\pi_*(X,\omega_{X/Y})$.
By (\ref{da0750}), there exists $\xi\in R^0\pi_*(X,{\cal
O}_X(nE_P))\otimes {\cal O}_Y(-C_0+C)$, such that
\begin{equation}\label{da0740}\xi=
\zeta_{n_0,n_0}-w_2
\end{equation} in $\pi_*({\cal O}_X(nE_P)/{\cal O}_X)\otimes F$.
 So our Theorem is true.
\end{proof}

\hspace{1in}

Let $m_2,m_3$ be positive integers. Let $\Omega_{X/Y}$ be the
sheaf of relative differentials of $X$ over $Y$. Let $\omega$ be a
section of $\Omega_{X/Y}$ on $X$, such that there exists sections
$u_i$ of ${\cal O}_X(m_2E_P)$ on $X$ and sections $e_i$ of ${\cal
O}_X(m_3E_P)$ on $X$, where $i\in I_1$, satisfying:
\begin{equation}\label{da6000}\sum_{i\in I_1}u_i\cdot e_i=0
\end{equation} and
\begin{equation}\label{da6010}\sum_{i\in I_1}\frac{u_i}{\eta^{m_2}}
\cdot d\frac{e_i}{\eta^{m_3}}=\omega
\end{equation} over $X$.

Let $n_0>9g^2$ be an integer. Let $N_0>0$ be the constant
determined by $n_0,X_{\sca{C}}$ constructed in {\bf Theorem
\ref{dh0800}}. We assume (\ref{ts4770}) is true in the rest of
this section. Let $\xi_0$ be a section of ${\cal
O}_{X}(n_0E_P)\otimes {\cal O}_Y(-C_0+C)$ that satisfies the
properties in {\bf Theorem \ref{dh0800}}, and let $C_0$ be the
divisor determined by $\xi_0=0$ as a section of ${\cal
O}_X(n_0E_P)$ on $E_P$. For simplicity, we assume $C_0=0$ and
$C=0$ in the rest of this section. The general case can be proved
with small modifications.

\begin{lemma}\label{da6100}
For $n>9g^2$, there exists a set of generators $\{e_{j,1}\}$ of
\begin{equation}\label{da4945}
R^0\pi_*(X, {\cal O}_X(nE_P))\end{equation} over $\Bbb{Z}$, such
that
\begin{equation}\label{da4950}\log\|e_{j,1}\|<\frac{n+1}{
[F:\ca{Q}]}\cdot \omega_{X/Y}\cdot E_P+O(1)
\end{equation}
\end{lemma}

\begin{proof}
By the assumption $C=0$ and  the arguments in {\bf Theorem
\ref{dh0800}}, the natural map by restriction
\begin{equation}\label{da6120}
R^0\pi_*(X,{\cal O}_X(n E_P))\longrightarrow \pi_*({\cal O}_X(n
E_P)/{\cal O}_X(3g E_P))
\end{equation} is surjective. Note there exists
a set of elements $\{e_{j,2}\}$ in
\begin{equation}\label{da6130}
\pi_*({\cal O}_X(n E_P)/{\cal O}_X(3g E_P))
\end{equation} that generates (\ref{da6130}) over $\ca{Z}$, and
satisfies \begin{equation}\label{da6140} \log\|e_{j,2}\|<\frac{n}{
[F:\ca{Q}]}\cdot \omega_{X/Y}\cdot E_P+d_Y+O(1)\end{equation} Then
by the existence of $\{e_i'\}$ that generates $R^0\pi_*(X,{\cal
O}_X(3g E_P))$ satisfying (\ref{da0840}), where $C_0$ is taken to
be zero, there exists a set of elements $\{e_{j,3}\}$ in
(\ref{da4945}), that generates (\ref{da4945}) over $\ca{Z}$, and
satisfies
\begin{equation}\label{da6160}
\log\|e_{j,3}\|<\frac{n}{ [F:\ca{Q}]}\cdot \omega_{X/Y}\cdot
E_P+d_Y+O(1)\end{equation} By the assumption that
\begin{equation}d_Y<\frac{1}{ [F:\ca{Q}]}\cdot \omega_{X/Y}
\cdot E_P-N_0 \label{da6162}\end{equation} we see our Lemma is
true.
\end{proof}

\begin{lemma}\label{da6999}
Let $N_1$ be the positive number satisfying
\begin{equation}\label{da6998}\log N_1=\frac{1}{[F:\ca{Q}]}
\cdot\omega_{X/Y}\cdot E_P
\end{equation}
 Let $\{\varrho_i:i\in I\}$ be the set of
all the complex numbers satisfying
\begin{equation}\label{da6980}
\varrho_i^{2n_0}=1\end{equation}

Then there exists sections $\xi_1,\xi_2$ of ${\cal
O}_X(n_1E_P),{\cal O}_X(n_2E_P)$ over $X$ respectively, and a set
of open sets $\{U_{i,\sigma}:i\in I,\sigma\in S\}$ on
$X_{\sca{C}}$, that satisfies the following:
\begin{enumerate}\setcounter{rom}{0}\item[\rmn] $n_1=2n_0$ and
 $\xi_1-\xi_0^{2}$ vanishes along $E_P$ as a section
of ${\cal O}_X(n_1E_P)$ on $X$. \item[\rmn] $n_2=3n_0$ and
$\xi_2-\xi_0^{3}$ vanishes along $E_P$ as a section of ${\cal
O}_X(n_2E_P)$ on $X$. \item[\rmn] Let $r_1>0$ be the number
constructed for $\{\varrho_i:i\in I\}$ in {\bf
Lemma~\ref{da3098}}. Then there exists a positive number $\kappa$
and a simple connected open set $U_{i,\sigma}\subset X_\sigma$ for
all $i\in I$ and $\sigma\in S$, such that
\begin{equation}\label{da6970} \left|\tau(x)-\varrho_i\right|
< r_1\end{equation} for all $x\in U_{i,\sigma}$, where
\begin{equation}\label{da6997} \tau=\frac{2\cdot
\xi_2}{\kappa^\frac{3}{2}\cdot
N_1^{n_2}\cdot\eta^{n_2}}\end{equation} on $X_{\sca{C}}$, and
\begin{equation}\label{da6960}\left|z(x)\right|\geq \kappa
\end{equation} for all $x\in X_{\sca{C}}\smallsetminus
\bigcup_{\sigma\in S}\bigcup_{i\in I}U_{i,\sigma}$, where
\begin{equation}\label{da8101} z=\frac{\xi_1}{
N_1^{n_1}\cdot\eta^{n_1}}\end{equation} on $X_{\sca{C}}$.
\item[\rmn] For $\sigma\in S$, let $\{x_{i,\sigma}:i\in I\}$ be
the set of all the  points on $X_\sigma$ where $\xi_1=0$. Then for
all $\sigma\in S$ and $i\in I$, we have $x_{i,\sigma}\in
U_{i,\sigma}$ and
\begin{equation}\label{da6950}\log
\left|\tau(x_{i,\sigma})-\varrho_i\right| <-\log
N_1+O(1)\end{equation} \item[\rmn] $dz(x)\neq 0$ for all $x\in
\bigcup_{\sigma\in S}\bigcup_{i\in I}U_{i,\sigma}$.
\end{enumerate}
\end{lemma}

\begin{proof}
Let $n_1=2n_0$ and $n_2={3}n_0$.  Let
$\{s_{1},\cdots,s_{n_1-g-1}\}$ be a set of elements in
\begin{equation}\label{da7090}
R^0\pi_*(X,{\cal O}_X((n_1-2)E_P))\otimes F\end{equation} that
generates (\ref{da7090}) over $F$. Let
$\{s_1',\cdots,s'_{n_2-g-1}\}$ be a set of elements in
\begin{equation}\label{da7130}
R^0\pi_*\left(X,{\cal
O}_X\left(\left(n_2-2\right)E_P\right)\right)\otimes
F\end{equation} that generates (\ref{da7130}) over $F$. Let
$\ca{P}^{n_1-g}_{\sca{C}}$ be the $n_1-g$ dimensional projective
space associated to $R^0\pi_*(X,{\cal O}_X(n_1E_P))\otimes_\sigma
\ca{C}$. Let $\ca{P}^{n_2-g}_{\sca{C}}$ be the $n_2-g$ dimensional
projective space associated to $R^0\pi_*(X,{\cal
O}_X(n_2E_P))\otimes_\sigma \ca{C}$.

Let $X_3\subset \ca{P}^{n_1-g}_{\sca{C}}\times
\ca{P}^{n_2-g}_{\sca{C}}$ be the algebraic closure of the set of
points
\begin{equation}\label{da7140}
\left(\xi_0^2+\sum_{i=1}^{n_1-g-1}z_i\cdot s_i, \ \
\xi_0^3+\sum_{i=1}^{n_2-g-1}z_i'\cdot s_i'\right)\end{equation}
where $ z_i,z_i'\in\ca{C}$ and
$\xi_0^2+\sum_{i=1}^{n_1-g-1}z_i\cdot s_i$ does not vanish at any
point on $X_\sigma$ with order greater than $1$, that satisfies
the following relation:

Let $\{x_{i,\sigma}':i\in I\}$   be the set of all the points in
$X_\sigma$ satisfying
\begin{equation}\label{da7150}
\xi_0^2(x_{i,\sigma}')+\sum_{j=1}^{n_1-g-1}z_j\cdot s_j
(x_{i,\sigma}')=0\end{equation} Then
$\xi_0^3+\sum_{j=1}^{n_2-g-1}z_j'\cdot
s_j'-\frac{\kappa^{\frac{3}{2}}\cdot\varrho_i \cdot
N_1^{n_2}}{2}\cdot \eta^{n_2}$ vanishes at point $x_{i,\sigma}'$
with order $2$ for all $\sigma\in S$ and $i\in I$.

For $\xi_0^2+\sum_{i=1}^{n_1-g-1}z_i\cdot s_i$, let
$x_{i,\sigma}'$ be a point in $X_\sigma$ where
$\xi_0^2+\sum_{i=1}^{n_1-g-1}z_i\cdot s_i$ is equal to zero. Then
to find an element $\xi_0^3+\sum_{i=1}^{n_2-g-1}z_i'\cdot s_i'$
such that
\begin{equation}\label{da7160}
\xi_0^3+\sum_{i=1}^{n_2-g-1}z_i'\cdot
s_i'-\frac{\kappa^{\frac{3}{2}}\cdot \varrho_i\cdot
N_1^{n_2}}{2}\cdot \eta^{n_2}
\end{equation} vanishes at point $x_{i,\sigma}'$ with order
$2$, we need two algebraic relations. So $X_3$ is of codimension
$2n_1+2$ in $\ca{P}^{n_1-g}_{\sca{C}}\times
\ca{P}^{n_2-g}_{\sca{C}}$. So $X_3$ is an algebraic variety of
dimension $n_0-2g-2>0$. Hence there exists a point
$(\xi_{5,\sigma},\xi_{6,\sigma})$ on $X_3$ that satisfies the
following:
\begin{enumerate}\item $\xi_{5,\sigma}$ is equal to $\xi_0^2$ at
point $E_P\otimes_\sigma\ca{C}$ as a section of ${\cal
O}_X(n_1E_P)$ on $X_\sigma$. And
\begin{equation}\label{da7170}
d\frac{\xi_{5,\sigma}}{\eta^{n_1}}(x)\neq 0\end{equation} for all
the $x\in X_\sigma$ satisfying \begin{equation}\label{da7180}
\log\left|\frac{\xi_{5,\sigma}}{\eta^{n_1}}(x)\right|
<n_1\cdot\log N_1+O(1)\end{equation} \item $\xi_{6,\sigma}$ is
equal to $\xi_0^3$ at point $E_P\otimes_\sigma\ca{C}$ as a section
of ${\cal O}_X(n_2E_P)$ on $X_\sigma$.
\end{enumerate}

Let $\{x_{i,\sigma}'':i\in I\}$ be the set of all the points on
$X_\sigma$, where $\xi_{5,\sigma}=0$. Since $\xi_{6,\sigma}-
\frac{1}{2}\cdot \kappa^{\frac{3}{2}}\cdot \varrho_i\cdot
N_1^{n_2}\cdot\eta^{n_2}$ vanishes at $x_{i,\sigma}''$ with order
$2$, so
\begin{equation}\label{da7190}
\frac{\xi_{6,\sigma}}{N_1^{n_2}\cdot
\eta^{n_2}}=\frac{\kappa^{\frac{3}{2}}\cdot
\varrho_i}{2}+O(\kappa^2)
\end{equation} over the simple connected open set containing
point $x_{i,\sigma}''$ determined by
\begin{equation}\label{da7200}\left|
\frac{\xi_{5,\sigma}}{N_1^{n_1}\cdot\eta^{n_1}}\right|<2\kappa
\end{equation}
when $\kappa$ is small enough. For all $\sigma\in
S$, choose $\xi_{5,\sigma},\xi_{6,\sigma}$ suitably, such that
$\xi_{5,\sigma},\xi_{6,\sigma}$ are mapped to
$\xi_{5,\overline{\sigma}},\xi_{6,\overline{\sigma}}$ respectively
under the complex conjugation, where $\overline{\sigma}$ denotes
the complex conjugate of $\sigma$. Then there exists elements
\begin{equation}\label{da7210}
\xi_5\in R^0\pi_*(X,{\cal O}_X(n_1E_P))\otimes\ca{R}\end{equation}
\begin{equation}\label{da7220}
\xi_6\in R^0\pi_*(X,{\cal O}_X(n_2E_P))\otimes\ca{R}\end{equation}
such that $\xi_5=\xi_{5,\sigma}$ and $\xi_6=\xi_{6,\sigma}$ for
all $\sigma\in S$.

By {\bf Lemma \ref{da6100}}, there exists
\begin{equation}\label{da7230}
\xi_1\in R^0\pi_*(X,{\cal O}_X(n_1E_P))\end{equation}
\begin{equation}\label{da7240}
\xi_2\in R^0\pi_*(X,{\cal O}_X(n_2E_P))\end{equation}  such that
$\xi_1=\xi_0^2$ and $\xi_2=\xi_0^3$ along $E_P$ as sections of
${\cal O}_X(n_1E_P)$ and ${\cal O}_X(n_2E_P)$ on $X$ respectively,
and
\begin{equation}\label{da7250}\log\|\xi_1-\xi_5\|
<(n_1-1)\log N_1+O(1)\end{equation}
\begin{equation}\label{da7260}
\log\|\xi_2-\xi_6\|< (n_2-1)\log N_1+O(1)
\end{equation}
For $\sigma\in S$, let $\{x_{i,\sigma}:i\in I\} $ be the set of
points on $X_\sigma$, where $\xi_1=0$.

By (\ref{da7250}) (\ref{da7260}), we have
\begin{equation}\label{da7270}
\log\frac{\|\xi_1-\xi_5\|}{\|\xi_1\|} <-\log
N_1+O(1)\end{equation} \begin{equation}\label{da7280}
\log\frac{\|\xi_2-\xi_6\|}{\|\xi_2\|} <-\log
N_1+O(1)\end{equation} Since $-\log N_1\ll 0$, so by
(\ref{da7190}) (\ref{da7200}), we have
\begin{equation}\label{da7290}
\frac{2\cdot\xi_{2}}{\kappa^{\frac{3}{2}}\cdot
N_1^{n_2}\cdot\eta^{n_2}}=
 \varrho_i+O\left(1\right)\cdot \kappa^{\frac{1}{2}}
\end{equation} over the simple connected open set containing
point $x_{i,\sigma}$ determined by
\begin{equation}\label{da7300}\left|
\frac{\xi_{1}}{N_1^{n_1}\cdot\eta^{n_1}}\right|<\kappa
\end{equation} for all $i\in I$ and $\sigma\in S$,
when $\kappa$ is small enough. Let $U_{i,\sigma}$ be the simple
connected open set on $X_\sigma$ containing $x_{i,\sigma}$
determined by (\ref{da7300}).

(\trm{1}) (\trm{2}) are true by the definitions of $\xi_1,\xi_2$.
(\trm{3}) is implied by (\ref{da7290}) (\ref{da7300}) when
$\kappa$ is small enough. (\trm{4}) is implied by (\ref{da7280}).
(\trm{5}) is implied by (\ref{da7170}) (\ref{da7180}) when
$\kappa$ is small enough. So our Lemma is true.
\end{proof}

\hspace{1in}

Let $n_1,n_2,\xi_1,\xi_2,N_1,r_1,\kappa$ be the elements defined
in {\bf Lemma \ref{da6999}}. Let $\ca{P}_Y^1$ be the projective
line over $Y$. Let $\varphi_3:X_{\sca{C}}\longrightarrow
\ca{P}_Y^1\otimes\ca{C}$ be the morphism defined by
$[\eta^{n_2}:\xi_2]$. Let $v_0$ be the section of the canonical
line bundle ${\cal O}_{\sca{P}_Y^1}(1)\otimes\ca{C}$ that is equal
to $\eta^{n_2}$. Let $v_1$ be the section of the canonical line
bundle ${\cal O}_{\sca{P}_Y^1}(1)\otimes\ca{C}$ that is equal to
$\frac{2\cdot\xi_2}{\kappa^{\frac{3}{2}}\cdot N_1^{n_2}}$. For
convenience we still use $v_0,v_1$ to denote their pull backs on
$X_{\sca{C}}$ under $\varphi_3$ respectively.

Let $\tau,z$ be the elements defined in {\bf Lemma \ref{da6999}
(\trm{3})}. Let $D_1$ be the divisor determined by $\xi_1=0$ on
$X$. Assume $D_1$ is a horizontal divisor on $X$. Since $D_1$ does
not intersect with $E_P$ on $X$, so there exists a morphism
\begin{equation}\varphi_4:X\longrightarrow\ca{P}_Y^1
\label{da8129}\end{equation} defined by $[\xi_1:\eta^{n_1}]$. Let
$b_{m,l}$ be the numbers defined in (\ref{da8910}). Let $G$ be the
rational function on $X_{\sca{C}}$ defined by
\begin{equation}G= \sum_{i\in I_1}\sum_{l=0}^{m}b_{m,l}\cdot
\frac{v_1^{l}}{v_0^{l}}\cdot \frac{u_i}{\eta^{m_2}}\cdot
\varphi_4^*\left(\frac{\partial}{\partial z}{\rm
Trace}_{\varphi_4}\left( \frac{v_1^{m-l}}{v_0^{m-l}} \cdot
\frac{e_i}{\eta^{m_3}}\right) \right)\label{da8130}\end{equation}

\begin{theorem}\label{da2500}
Let $I,U_{i,\sigma}$ be the elements in {\bf Lemma \ref{da6999}}.
For $x\in X_{\sca{C}}\smallsetminus E_P\otimes\ca{C}$ satisfying
$d\frac{\xi_1}{\eta^{n_1}}(x)\neq 0$, assume
\begin{equation}\label{da8240} \omega-\beta_x\cdot
d\frac{\xi_1}{\eta^{n_1}}\end{equation} vanishes at point $x$,
where $\beta_x\in\ca{C}$.  Then there exists $\rho\in (0,1)$ and
$N_4>0$, such that for all $m>N_4$ and $x\in \bigcup_{\sigma\in
S}\bigcup_{i\in I}U_{i,\sigma}$, we have
\begin{equation}\label{da8099}
\left|G(x)-\frac{v_1^m}{v_0^m}(x)\cdot{\beta_x}\right|<
\left|\frac{v_1}{v_0}(x)\right|^m\cdot\rho^m
\end{equation}
\end{theorem}

\begin{proof}
This is implied by {\bf Theorem \ref{da9000}} and {\bf Lemma
\ref{da6999} (\trm{3})}.
\end{proof}

\begin{theorem}\label{da4912}
Let $p$ be the closed point in $X\otimes F$ determined by
$E_P\otimes F$. Assume there exists a set of elements
$\{e_{p,0},e_{p,1},\cdots\}$ in $R^0\pi_*(X, {\cal
O}_X(n_1E_P))\otimes F$ that satisfies the following:
\begin{enumerate} \item There exists elements $\{u_{p,0},u_{p,1},
\cdots\}$ in $R^0\pi_*(X,{\cal O}_X(m_2E_P))\otimes F$, such that
\begin{equation}\label{da8520}\sum_{i\in I_1}u_i\otimes e_i=
\sum_iu_{p,i}\otimes e_{p,i}
\end{equation} over $F$.
\item $e_{p,0}-\xi_1$ vanishes at point $p$  as a section of
${\cal O}_X(n_1E_P)\otimes F$ on $X\otimes F$, and
$\frac{u_{p,0}}{\eta^{m_2}}=1$. \item $e_{p,i}$ vanishes at point
$p$ as a section of ${\cal O}_X(n_1E_P)\otimes F$ on $X\otimes F$
for all $i>0$. \item $m_2<n_0$.
\end{enumerate}

Let $m$ be a positive even integer. Then there exists a section
$G_1$ of ${\cal O}_X(3n_0mE_P)$ on $X$, such that
\begin{equation}\label{da8511} \frac{G_1}{\eta^{{3n_0m}}}
= \frac{(2m)!\cdot N_1^{n_2m-n_1}\cdot \kappa^{\frac{3m}{2}}
}{m!\cdot m!\cdot 2^m}\cdot G\end{equation}  on $X_{\sca{C}}$, and
\begin{equation}\label{da8730}
\frac{G_1}{\eta^{3n_0m}}=N_2\cdot
\frac{\xi_1^{\frac{3m}{2}}}{\eta^{{3n_0m}}}+
\frac{t_1}{\eta^{{3n_0m}}}\end{equation} where $t_1$ is a section
of ${\cal O}_X({3n_0m}E_P)$  that vanishes along $E_P$ on $X$, and
$N_2$ is a positive integer satisfying
\begin{equation}\label{da8521}N_2>
\frac{(2m)!}{m!\cdot m!}
\end{equation}

Assume $\beta_x\neq 0$ for all the points $x$ on an open set
containing the closure of $\bigcup_{\sigma\in S}\bigcup_{i\in
I}U_{i,\sigma}$ in $X_{\sca{C}}$. Then there exists $N_3>0$, such
that
\begin{equation}\label{da8522}\left|\frac{\xi_1}{\eta^{n_1}}(x)
\right|\geq\kappa\cdot N_1^{n_1}
\end{equation} for all $x\in
X_{\sca{C}}$ where $G_1(x)=0$, when $m>N_3$.
\end{theorem}

\begin{proof} Note we have
\begin{equation}f_5\left(\frac{(2m)!}{m!\cdot m!}v_0^mv_1^m\right)
=\sum_{l=0}^m\frac{m!}{l!\cdot (m-l)!}\cdot \frac{m!}{l!\cdot
(m-l)!}\cdot v_0^{m-l}v_1^l\otimes
v_0^lv_1^{m-l}\label{da8500}\end{equation} So we have
\begin{eqnarray}\lefteqn{\frac{(2m)!\cdot N_1^{n_2m-n_1}\cdot
\kappa^{\frac{3m}{2}} }{2^m\cdot m!\cdot m!}\cdot G=\sum_{i\in
I_1}\sum_{l=0}^m \frac{m!}{l!\cdot (m-l)!}
}\nonumber\\
&& \cdot \frac{m!}{l!\cdot
(m-l)!}\cdot\frac{\xi_2^l}{\eta^{ln_2}}\cdot
\frac{u_i}{\eta^{m_2}} \cdot
\varphi_4^*\left(\frac{\partial}{\partial
\frac{\xi_1}{\eta^{n_1}}}{\rm Trace}_{\varphi_4}\left(
\frac{\xi_2^{m-l}}{\eta^{(m-l)n_2}} \cdot
\frac{e_i}{\eta^{n_1}}\right) \right)\hspace{.1in}
\label{da8510}\end{eqnarray} So (\ref{da8511}) is a section of
${\cal O}_X(n_5E_P)$ on $X$, where $n_5$ is a positive integer.

Let's consider the case that $l$ is an even number. Since
$\frac{\xi_2}{\eta^{n_2}}$ has a pole of order
$n_2=\frac{3n_1}{2}$ along $E_P$, so
$\frac{\xi_2^l}{\eta^{n_2l}}\cdot\frac{e_{p,i}}{\eta^{n_1}} $ has
a pole of order $< 3n_1\cdot\frac{l}{2}+n_1$  along $E_P$ for all
$i>0$. Therefore
\begin{equation}\label{da8560}{\rm Trace}_{\varphi_4}
\left(\frac{\xi_2^l}{\eta^{n_2l}}\cdot\frac{e_{p,i}}{\eta^{n_1}}
\right)=\sum_{j=0}^{\frac{3l}{2}}c_{j,l,i
}\cdot\frac{\xi_1^j}{\eta^{n_1j}}
\end{equation} for $i>0$, where $c_{j,l,i}\in F$. Since
$\xi_2=\xi_0^3$ and $\xi_1=\xi_0^{2}$ along $E_P$, so we have
\begin{equation}\label{da8570}
{\rm Trace}_{\varphi_4}
\left(\frac{\xi_2^l}{\eta^{n_2l}}\cdot\frac{e_{p,0}}{\eta^{n_1}}
\right)=\frac{n_1\cdot\xi_1^{\frac{3l}{2}+1}}{
\eta^{n_1(\frac{3l}{2}+1)}}
+\sum_{j=0}^{\frac{3l}{2}}c_{j,l,0}\cdot\frac{\xi_1^j}{\eta^{n_1j}}
\end{equation} where $c_{j,l,0}\in F$.

Now consider the case that $l$ is an odd number. Since
$\frac{\xi_2}{\eta^{n_2}}$ has a pole of order
$n_2=\frac{3n_1}{2}$ along $E_P$, so
$\frac{\xi_2^l}{\eta^{n_2l}}\cdot\frac{e_{p,i}}{\eta^{n_1}} $ has
a pole of order $\leq n_1\cdot\frac{3l+2}{2}$  along $E_P$ for all
$i$. Therefore
\begin{equation}\label{da8575}{\rm Trace}_{\varphi_4}
\left(\frac{\xi_2^l}{\eta^{n_2l}}\cdot\frac{e_{p,i}}{\eta^{n_1}}
\right)=\sum_{j=0}^{\frac{3l+1}{2}}
c_{j,l,i}'\cdot\frac{\xi_1^j}{\eta^{n_1j}}
\end{equation} where $c_{j,l,i}'\in F$.

Therefore by $m_2<\frac{n_1}{2}$ and (\ref{da8560}) (\ref{da8570})
and $\frac{u_{p,0}}{\eta^{m_2}}=1$, when $l$ is an even number, we
have
\begin{eqnarray}\lefteqn{\sum_i\frac{\xi_2^l}{\eta^{ln_2}}\cdot
\frac{u_{p,i}}{\eta^{m_2}} \cdot
\varphi_4^*\left(\frac{\partial}{\partial z}{\rm
Trace}_{\varphi_4}\left( \frac{\xi_2^{m-l}}{\eta^{(m-l)n_2}} \cdot
\frac{e_{p,i}}{\eta^{n_1}}\right) \right)}\nonumber\\
&=& \frac{n_1\cdot(\frac{3}{2}(m-l) +1)\cdot\xi_1^{\frac{3m}{2}}}{
\eta^{n_2m}} +G_{p,1}\hspace{1in}\label{da8580}\end{eqnarray}
where $G_{p,1}$ vanishes along $E_P\otimes F$ as a section of
${\cal O}_{X}(n_2mE_P)\otimes F$ on $X\otimes F$. By
$m_2<\frac{n_1}{2}$ and (\ref{da8575}), when $l$ is an odd number,
we have
\begin{equation}{\sum_i\frac{\xi_2^l}{\eta^{ln_2}}\cdot
\frac{u_{p,i}}{\eta^{m_2}} \cdot
\varphi_4^*\left(\frac{\partial}{\partial z}{\rm
Trace}_{\varphi_4}\left( \frac{\xi_2^{m-l}}{\eta^{(m-l)n_2}} \cdot
\frac{e_{p,i}}{\eta^{n_1}}\right)
\right)=0}\label{da8585}\end{equation} at point $E_P\otimes F$ as
a section of ${\cal O}_{X}(n_2mE_P)\otimes F$ on $X\otimes F$.

 By (\ref{da8510}) (\ref{da8580}) (\ref{da8585}), we have
\begin{equation}\label{da8600}
\frac{(2m)!\cdot N_1^{n_2m-n_1}\cdot \kappa^{\frac{3m}{2}}
}{2^m\cdot m!\cdot m!}\cdot G= N_2\cdot
\frac{\xi_1^{\frac{3m}{2}}}{\eta^{n_2m}}+\frac{G_2}{\eta^{n_2m}}
\end{equation} where  $G_2$ is a section of ${\cal
O}_X(n_2mE_P)$ vanishing along $E_P$ over $X$, and
\begin{equation}\label{da8620}N_2=n_1\cdot\sum_{l_1=0}^{
\frac{m}{2}}
\left(\frac{3}{2}(m-2l_1)+1\right)\cdot\frac{m!}{(2l_1)!\cdot
(m-2l_1)!}\cdot \frac{m!}{(2l_1)!\cdot (m-2l_1)!}
\end{equation}  Note
\begin{equation}\frac{(2m)!}{m!\cdot
m!}=\sum_{l=0}^{{m}} \frac{m!}{l!\cdot (m-l)!}\cdot
\frac{m!}{l!\cdot (m-l)!}\label{da8590}\end{equation} so we have
\begin{equation}\label{da8630}N_2>\frac{(2m)!}{m!\cdot m!}
\end{equation}
By {\bf Theorem~\ref{da2500}}, when $m\gg 0$, we have $G\neq 0$ on
$\bigcup_{\sigma\in S}\bigcup_{i\in I}U_{i,\sigma}$. Then by
(\ref{da6960}), we see (\ref{da8522}) is true. So our Theorem is
true.
\end{proof}

\begin{cor}\label{da2700}  Under all the conditions
assumed in {\bf Theorem \ref{da4912}}, we have
\begin{equation}\label{da2780}\frac{1}{n_1\cdot[F:\ca{Q}]}
\sum_{\sigma\in S}\sum_{j\in
I}\log\left|{G}(x_{j,\sigma})\right|\geq n_1\cdot\log N_1
+m\cdot\log 2
\end{equation}
\end{cor}

\begin{proof}
Let $D_1$ be the divisor determined by $\xi_1=0$ as a section of
${\cal O}_X(n_1E_P)$ on $X$. Let $D_3$ be the horizontal component
of the divisor determined by $G_1=0$ as a section of ${\cal
O}_X(n_2mE_P)$ on $X$. Let $C_2$ be the intersection cycle of
$D_1$ and $D_3$ on scheme $X$. Let $C_3$ be the the intersection
cycle of $E_P$ and $D_3$ on scheme $X$. Let $\{x_i:i\in I''\}$ be
the set of all the points on $X_{\sca{C}}$ where $G_1=0$. By {\bf
Theorem~\ref{da4912}}, we have
\begin{equation}\label{da2730}
\deg C_3\geq [F:\ca{Q}]\cdot\log\frac{(2m)!}{m!\cdot m!}
\end{equation}  Note $\{x_i:i\in I''\}$ has $n_2m\cdot[F:\ca{Q}]$
points. By (\ref{da8522}), we have
\begin{equation}\label{da2750}
\log\left|\frac{\xi_1}{\eta^{n_1}}(x_i)\right|> n_1\cdot\log N_1
+\log\kappa\end{equation} for all $i\in I''$. Therefore we have
\begin{eqnarray}\lefteqn{\frac{1}{[F:\ca{Q}]}\deg C_2=
\frac{1}{[F:\ca{Q}]}\sum_{i\in I''}
\log\left|\frac{\xi_1}{\eta^{n_1}}(x_i)\right|+\frac{n_1}{
[F:\ca{Q}]}\cdot\deg C_3}\nonumber\\
&\geq&n_2m\cdot\left({n_1}\log N_1+\log \kappa\right)
+n_1\cdot\log\frac{(2m)!}{m!\cdot m!}
\hspace{.8in}\label{da2760}\end{eqnarray} Since $E_P$ does not
intersect with $D_1$ on scheme $X$, so we have
\begin{equation}\label{da2770} \sum_{\sigma\in S}\sum_{j\in
I}\log\left|\frac{(2m)!\cdot N_1^{n_2m-n_1}\cdot
\kappa^{\frac{3m}{2}} }{2^m\cdot m!\cdot
m!}\cdot{G}(x_{j,\sigma})\right|=\deg C_2\end{equation} By
(\ref{da2770}) (\ref{da2760}), we see (\ref{da2780}) is true.
\end{proof}

\hspace{1in}

\begin{pot}{\ref{ts4200}} Note we have
\begin{equation}\label{da2710}d\left(
\frac{\xi_6^2}{\eta^{2n_6}}\cdot\frac{\xi_7}{\eta^{n_7}}\right)
-2\frac{\xi_6}{\eta^{n_6}}\cdot d\left(
\frac{\xi_6}{\eta^{n_6}}\cdot\frac{\xi_7}{\eta^{n_7}}\right)
+\frac{\xi_6^2}{\eta^{2n_6}}\cdot
d\left(\frac{\xi_7}{\eta^{n_7}}\right)=0 \end{equation} where
$\xi_6,\xi_7$ are sections of ${\cal O}_X(n_6E_P)$ and ${\cal
O}_X(n_7E_P)$ respectively. So we can replace $\sum_{i\in
I_1}u_i\otimes e_i$ by
\begin{equation}\label{da2720}\sum_{i\in
I_1}u_i\otimes e_i+\left(1\otimes \xi_6^2\xi_7-2\xi_6 \otimes
\xi_6\xi_7+\xi_6^2\otimes\xi_7\right)
\end{equation} By setting $n_7>2n_6$ and $n_0>m_2$ and by
choosing $\xi_6,\xi_7$, such that $\xi_6^2\xi_7$ is a section of
${\cal O}_X(n_1E_P)$ on $X$, and $\xi_6^2\xi_7-\xi_1=0$ along
$E_P$ as a section of ${\cal O}_X(n_1E_P)$ on $X$. Then
(\ref{da2720}) satisfies the four conditions required by
$\sum_{i\in I_1}u_i\otimes e_i$ in {\bf Theorem \ref{da4912}}. So
we assume $\sum_{i\in I_1}u_i\otimes e_i$ satisfies the four
conditions in {\bf Theorem \ref{da4912}}, and the canonical norm
of $\sum_{i\in I_1}u_i\otimes e_i$ is equal to $O(N_1^{n_1})$. By
choosing $\xi_1$ suitably, we assume $\beta_x\neq 0$ for all the
points $x$ on an open set containing the closure of
$\bigcup_{\sigma\in S}\bigcup_{i\in I}U_{i,\sigma}$ in
$X_{\sca{C}}$.

Let $m$ be an even integer satisfying $m>N_3,N_4$, where $N_3,N_4$
are defined in {\bf Theorem \ref{da4912} \ref{da2500}}. By
checking the constructions of $N_3,N_4$, we can assume
$N_3,N_4<N_1^{\frac{1}{2}}$. Let $m$ be an even integer in
$(N^{\frac{1}{2}},2N_1^{\frac{1}{2}})$. By (\ref{da8240})
(\ref{da8099}) (\ref{da6950}), we have
\begin{equation}\label{da2800} \log\|\omega\|\geq
\frac{1}{n_1\cdot [F:\ca{Q}]}\sum_{\sigma\in S}\sum_{j\in
I}\log\left|{G}(x_{j,\sigma})\right|+O(1)\end{equation} where
$\|\omega\|$ denotes the canonical norm of $\omega$. By
(\ref{da2780}), we have
\begin{equation}\label{da2810}\log\|\omega\|\geq n_1\cdot
\log N_1 +N_1^{\frac{1}{2}}\cdot\log 2+O(1)
\end{equation} for a non zero section $\omega$ of $\Omega_{X/Y}$
on $X$ that is equal to  $\sum_i\frac{u_i}{\eta^{m_2}}\cdot
d\frac{e_i}{\eta^{m_3}}$. Let $n_1\longrightarrow +\infty$, we see
$\|\omega\|\longrightarrow+\infty$. This is not possible. So there
exists a constant $N_0>0$, determined by $X_{\sca{C}}$, such that
(\ref{ts4770}) is not true. Therefore {\bf Theorem~\ref{ts4200}}
is true. \end{pot}

\end{document}